\documentclass{ifacconf}

\usepackage{graphicx}      
\usepackage{natbib}        

\usepackage{amsfonts,amssymb,amsmath}
\usepackage[update,prepend]{epstopdf}

\newcommand{\beq}{\begin{equation}}
\newcommand{\beqnt}{\begin{equation}\nonumber}
\newcommand{\eeq}{\end{equation}}
\newcommand{\fref}[1]{{\rm(\ref{#1})}}                      
\newcommand{\PC}{{\cal P}}
\newcommand{\QC}{{\cal Q}}
\newcommand{\UA}{{\mathbb U}}
\newcommand{\UB}{{\mathbf U}}
\newcommand{\VB}{{\mathbf V}}
\newcommand{\Upr}{{\cal U}}
\newcommand{\Vpr}{{\cal V}}
\newcommand{\Uqs}{{\mathbf Q}}
\newcommand{\Uxu}{{\mathbf S}}
\newcommand{\Part}{\Delta}
\newcommand{\Partk}{{\Part_k}}
\newcommand{\Pparset}{{\Part_T}}
\newcommand{\nPart}{n_{\Part}}
\newcommand{\cPart}{(\tau_i)_{i\in\nint{0}{\nPart}}}
\newcommand{\Ul}{\UA_*}
\newcommand{\Ulpar}{{\UB_*^{\Part}}}       
\newcommand{\Ulp}[1]{{\UB^{\Part}_{*#1}}}   
\newcommand{\epar}{\varepsilon}
\newcommand{\Ue}{{\UA_\epar}}
\newcommand{\Uepar}{{\UB^{\Part}_\epar}}       
\newcommand{\Uep}[1]{{\UB^{\Part}_{\epar#1}}}   
\newcommand{\RA}{{\mathbb R}}
\newcommand{\Lp}{{\ensuremath{L_{\mathrm p}}}}  

\newcommand{\Xc}{{\cal X_\textsc{c}}}                
\newcommand{\Xs}{{\cal X}}                           
\newcommand{\Xpr}{{\cal X_\textsc{p}}}               
\newcommand{\Gs}{{\Gamma}}                           
\newcommand{\Gc}{{\Gamma_\textsc{c}}}                         
\newcommand{\Gp}{{\Gamma_\textsc{p}}}                         
\newcommand{\Gq}{{\Gamma_\textsc{q}}}                         
\newcommand{\Wq}{{\cal W}}                         
\newcommand{\qndx}{\gamma}                          
\newcommand\argmin{\operatornamewithlimits{\mathrm{argmin}}}
\newcommand\argmax{\operatornamewithlimits{\mathrm{argmax}}}
\newcommand\diam[1]{\ensuremath{\operatorname{{\mathrm D}(#1)}}}

\newcommand\close{\operatorname{\mathbf cl}}

\newcommand*{\LL}[1]{\left#1}
\newcommand{\RR}[1]{\right#1}
\newcommand\comp[2]{{\mathbf{comp}_{#2}{(#1)}}}

\newcommand{\SP}[2]{{\left\langle{#1},{#2}\right\rangle}}
\newcommand{\nint}[2]{{{#1}..{#2}}}
\newcommand{\mydef}{=}
\newcommand\mysim[1]{\mathrel{\mathop{\sim}\limits_{#1}}}
\newcommand{\myinf}{\operatornamewithlimits{\inf\vphantom{\sup}}}
\newcommand{\mymin}{\operatornamewithlimits{\min\vphantom{\sup}}}
\newcommand{\mymax}{\operatornamewithlimits{\max\vphantom{\sup}}}

\begin{document}
\begin{frontmatter}

\title{Guaranteed control design under $\Lp$-compact  constraints on the disturbance} 

\author[First]{Dmitrii A. Serkov} 

\address[First]{Krasovskii Institute of Mathematics and Mechanics \\
Ural Branch, Russian Academy of Sciences\\
S.Kovalevskoy, 16, Yekaterinburg, 620990, Russia\\
and\\
Yeltsyn Ural Federal University\\
Mira, 19, Yekaterinburg, 620000, Russia}

\begin{abstract}                
The paper deals with the problem of optimization of a guaranteed (worst case) result for a control system described by an ordinary differential equation. 
The disturbances as functions of time are subject to functional constraints belonging to a given family of constraints. The latter family is known to the side that forms the control actions. The controlling side uses positional full-memory strategies and does not observe the disturbance. 
When the constraints family consists of $\Lp$-compact sets the optimal guaranteed result is non-improvable in the sense that it coincides with that obtained in the class of quasi-strategies -- nonanticipatory transformations of disturbances into controls. 

In this paper for the effectiveness of implemented control algorithm an additional condition on the system and appropriate ways of constructing an optimal strategy are specified.
\end{abstract}


\end{frontmatter}

\section{Introduction}
This work is related to Krasovskii's theory of guaranteeing positional control (see \cite{KraSub88e}, \cite{SubChe81e}). The theory focuses on assessment of the optimal guaranteed result -- the minimax of a cost functional -- for the control side that opposes the disturbance in the process of steering a dynamical control system. The properties of the optimal guaranteed result and an optimal strategy for the control side in the case where the non-observable dynamical disturbance subject to a functional constraint belonging to a given family of constraints.

Control problems with additional functional constraints imposed on the input dynamical disturbances have numerous interpretations and have been studied in various formalizations. The simplest functional constraint restricts the disturbances to the open-loop ones. In \cite{KrasIZVD70e}, \cite{Kra71DANe}, \cite{KrasUDS85e} the maximin open loop constructions (including the stochastic ones) use open-loop disturbances to find the optimal guaranteed result and optimal closed-loop strategies in control problems with non-constrained disturbances. 
In \cite{BarSub70e}, \cite{BarSub71e} properties of linear control systems in the cases of open-loop disturbances, disturbances generated by continuous feedbacks, and disturbances formed by upper semicontinuous set-valued closed-loop strategies were compared. 

In \cite{Kry91} assuming that the disturbances are restricted to an unknown $L_2$-compact set, it was shown that the optimal guaranteed result achieved in the class of the full-memory closed-loop control strategies equals that achieved by the 'fully informed controller' allowed to control the system using quasi-strategies --- nonanticipatory open-loop control responses to disturbance realizations \cite{SubChe81e}; in this sense the full-memory closed-loop strategies are uninprovable. In \cite{Ser_DAN2013et}, considering the problem setting proposed in \cite{Kry91} in the case of a continuous cost functional, new unimprovability conditions for full-memory control strategies were given and an optimal full-memory control strategy allowing numerical implementation was constructed. 

For the case of a continuous cost functional in \cite{Ser_TRIMM2014e}   it was shown that, firstly, the guaranteed control problem with open-loop disturbances is equivalent to that with the disturbances restricted to $L_2$-compact sets, and, secondly, the optimal guaranteed results achieved by the controlling player in the class of full-memory control strategies under these two types of constraints on the disturbances are equal to that achieved in the class of quasi-strategies. In showing the results, the elements of theory of robust dynamical inversion of control systems (\cite{KryOsi83e}, \cite{KryOsi95}) were used. However,  the control algorithm used in that study is not suitable for numerical implementation.

In this paper a new weakened condition on the control system (as compared with that in \cite{Kry91} and \cite{Ser_DAN2013et}), allowing numerical realization of an optimal full-memory control strategy, is provided.

Plan of the paper is as follows: in section 2 we give a formal statement of the problem and note the conditions of the works \cite {Kry91} and \cite{Ser_DAN2013et}, in section 3, we present a solution of the problem, corresponding to the work \cite {Ser_TRIMM2014e}, and section 4 shall give the weakened conditions on the system and a modernization of the solutions suitable for further numerical implementation.

\section{Definitions}

Consider a control system
\begin{equation}\label{sys}
\begin{cases}
\dot x(\tau)=f(\tau,x(\tau),u(\tau),v(\tau)),&\tau\in T\mydef[t_0,\vartheta]\subset\RA,\\
x(t_0)=z_0\in G_0\subset\RA^n,
\end{cases}
\end{equation}
$$
u(\tau)\in\PC\subset\RA^p,\ v(\tau)\in\QC\subset\RA^q,\quad\tau\in T.
$$
Here
$\PC$, $\QC$, and $G_0$ are compact sets;
and
$f(\cdot): T \times \RA^n \times \PC \times \QC \mapsto \RA^n $ is continuous,
locally Lipschitz in the second argument and such that for some $ K\ge0$ the inequality
$$
\sup_{(\tau,u,v)\in T\times\PC\times\QC}\|f(\tau,x,u,v)\|\le K(1+\|x\|)
$$
holds for all $x\in\RA^n$ ($\|\cdot\|$ denotes the norm in an Euclidian space).
{\em Controls} $u(\cdot): T \mapsto \PC $ and {\em disturbances} $v(\cdot): T \mapsto \QC $ are supposed to be  Lebesgue measurable.
Denote by $\Upr$ the set of all controls and by $\Vpr$ the set of all disturbances.

For arbitrary $(t_*,x_*)\in T\times\RA^n$, $u(\cdot)\in\Upr$, $v(\cdot)\in\Vpr$ we denote
by $x(\cdot,t_*,z_*,u(\cdot),v(\cdot))$
the (unique) Carath\'eodory solution of (\ref{sys}) (see \cite[II.4]{Wargae}) defined on $[t_*,\vartheta]$ and satisfying the initial condition $x(t_*)=x_*$.
We fix a compact set $G\subset T\times\RA^n$
such that
$ (t,x(t,t_0,z_0,u(\cdot),v(\cdot))) \in G $
for all
$ t \in T $,
$ z_0 \in G_0 $
$u(\cdot)\in\Upr$,
$v(\cdot)\in\Vpr$.

A set $\Part\mydef(\tau_i)_{i\in\nint{0}{\nPart}}$
where
$\tau_0=t_0$, $\tau_{i-1}<\tau_{i}$, $\tau_{\nPart}=\vartheta$
will be called a {\em partition} (of
interval $T$).
Denote by $\Pparset$ the set of all partitions.
For a partition $\Part\mydef(\tau_i)_{i\in\nint{0}{\nPart}}$ and
a $t\in T$ set $i_t\mydef\max_{i\in\nint{0}{\nPart}\atop\tau_i\le t}i$,
$\diam{\Part}\mydef\max_{i\in\nint{1}{\nPart}}\tau_i-\tau_{i-1}$.

Following \cite{Kry91}, define full-memory control strategies
used by the controlling player.
For every
$ \tau_*, \tau^* \in T $
where
$ \tau^* > \tau_* $
denote
by
$ {\cal U}|_{[\tau_*, \tau^*)} $
the set of the restrictions of
all controls to
$ [\tau_*, \tau^*) $.
Given a partition  $\Part\mydef(\tau_i)_{i\in\nint{0}{\nPart}}$,
any family
$\UB^\Part\mydef(\UB^\Part_i(\cdot))_{i\in\nint{0}{(\nPart-1)}}$,
where
$
\UB^\Part_i(\cdot): C([t_0, \tau_i], \RA^n) \mapsto {\cal U}|_{[\tau_i, \tau_{i+1})} $
$ (i\in\nint{0}{(\nPart-1)}) $,
will be called a
\emph{full-memory feedback for partition $\Part$}.
Every family
$ \UA = (\UB^{\Part})_{\Part \in \Pparset} $
where
$ \UB^{\Part} $
is a
full-memory feedback for $\Part$
will be called a
{\em full-memory control strategy}
(for the controlling player).
We denote by
$ \Uxu $
the set of all full-memory control strategies.

Given
a $ z_0 \in G_0 $,
a partition
$\Part\mydef(\tau_i)_{i\in\nint{0}{\nPart}}$,
a full-memory feedback
$ \UB^\Part = (\UB^\Part_i (\cdot))_{i \in \nint{0}{(n_{\Part}-1)}} $
for
$\Part $
and
a disturbance $ v(\cdot) \in \Vpr $,
the
function
$ x(\cdot) = x(\cdot, t_0, z_0, u(\cdot), v(\cdot)) $
where
$ u(\cdot) \in \Upr $
is such that
$ u(t) = $
$ \UB^\Part_{i_t}(x(\cdot)|_{[t_0,\tau_{i_t}]})(t) $
for all
$ t \in T$
will be called
the (system's) {\em motion} originating at $ z_0 $ and corresponding to $\Part$,  $ \UB^\Part $ and $ v(\cdot) $;
we denote $ x(\cdot) $ and $ u(\cdot) $ by $ x(\cdot,z_{0},\UB^{\Part},v(\cdot)) $ and  $ u(\cdot,z_{0},\UB^{\Part},v(\cdot)) $, respectively.

For every
$ z_0 \in G_0 $,
every full-memory control strategy
$ \UA = (\UB^{\Part})_{\Part \in \Pparset} \in \Uxu $
and every nonempty set
of disturbances,
$\VB \subseteq \Vpr $,
we define the
{\em bundle of
motions}
originating from
$ z_0 $
and corresponding to
$ \UA $
and
$\VB $
to be the set 
$ X(z_0,\UA,\VB) $
of all
$ x (\cdot) \in C (T; \RA^n) $
with the following property:
there is a sequence
$ \{(z_ {0k}, v_k (\cdot), \Partk, \UB^{\Partk}) \}_{k=1}^{\infty} $
in
$ G_0 \times \VB \times \Pparset \times \UA $
such that
$\lim_{k \to\infty} z_{0k} = z_0$, $\lim_{k \to\infty} \diam{\Partk} = 0$
and
$ x(\cdot,z_{0k},\UB^{\Partk},v_k(\cdot)) \rightarrow x(\cdot) $
in
$ C(T;\RA^n) $.

In the above definition,
$ \VB $
is a functional constraint on the disturbance. Generally, we assume that the controlling side does not know $ \VB $ but knows a class of functional constraints $ \VB $ belongs to. The latter class gives a general information on the disturbing constraints but does not provide information on the exact ones.

In the above cases,
for every
$ z_0 \in G_0 $
and every
full-memory control strategy
$ \UA \in \Uxu$,
we define the bundles of the
system's
{\em motions}
originating at
$ z_0 $
under
$ \UA $
subject to
{\em arbitrary disturbances},
{\em $\Lp$-compactly constrained disturbances} (for some fixed $\rm p>1$),
and
{\em open-loop disturbances}
as, respectively,
\begin{eqnarray*}
\Xs(z_0,\UA)&\mydef&X(z_0,\UA,\Vpr),\notag\\
\Xc(z_0,\UA)&\mydef&\bigcup_{\VB\in\comp{\Vpr}{\Lp(T,\RA^q)}}X(z_0,\UA,\VB),\\
\Xpr(z_0,\UA)&\mydef&\bigcup_{v(\cdot)\in\Vpr}X(z_0,\UA,\{v(\cdot)\}) ;\notag
\end{eqnarray*}
here $\comp{\Vpr}{\Lp(T,\RA^q)}$ denotes the family of all subsets of $ \Vpr$ compact in $\Lp(T; \RA^q)$.

\begin{rem}
The definition of $\Xs(z_0, \UA)$ is a straightforward generalization of the definition of the
set of constructive motions
generated by a closed-loop control strategy (see \cite{KraSub88e}).
The definition of $\Xc(z_0, \UA)$ follows \cite{Kry91}.
\end{rem}

According to the definitions, $\Xpr(z_0, \UA) \subseteq\Xc(z_0, \UA) \subseteq \Xs (z_0, \UA)$ holds for all $z_0\in G_0$ and $\UA\in\Uxu$.
In \cite {Ser_UDSU_09e} it was shown that generally $\Xpr(z_0, \UA)\neq\Xs (z_0, \UA)$.
A similar reasoning can lead to a statement that generally
$\Xc(z_0, \UA)\neq\Xs (z_0, \UA)$.

Let the controlling player evaluate the quality of the system's motions by a continuous
{\em cost functional}
$ \qndx(\cdot): C (T; \RA^n)\mapsto\RA $
acting as a benefit functional for the disturbing player.
The controlling player
seeks then to chose a full-memory control strategy that guarantees the minimum
value for the supremum of $ \qndx(x(\cdot)) $
over the system's motions
$ x(\cdot) $
corresponding to the chosen control strategy and all disturbances
that are allowed to be chosen by the disturbing player
within the given constraints.

Following \cite {KraSub88e}, and \cite {SubChe81e},
we call
$$
\Gs (z_0, \UA) \mydef \sup_ {x (\cdot) \in \Xs (z_0, \UA)} \qndx (x(\cdot))
$$
the {\em guaranteed result} at $ z_0 \in G_0 $ for a full-memory control strategy $\UA  $ against \emph{arbitrary disturbances}; and we call
$$
\Gs(z_0) \mydef \inf_{\UA \in \Uxu} \Gs (z_0, \UA)
$$
the {\em optimal guaranteed result} at $ z_0 \in G_0 $ in the class of the full-memory control strategies $\Uxu$, against \emph{arbitrary disturbances}.
Similarly, we call
$$
\Gc (z_0, \UA) \mydef \sup_ {x (\cdot) \in \Xc (z_0, \UA)} \qndx (x(\cdot))
$$
the {\em guaranteed result} at $ z_0 \in G_0 $ for a full-memory control strategy $\UA  $ against $\Lp$-\emph{compactly constrained disturbances} and we call
$$
\Gc(z_0) \mydef \inf_{\UA \in \Uxu} \Gc(z_0, \UA)
$$
the {\em optimal guaranteed result} at $ z_0 \in G_0 $ in $\Uxu$ against $\Lp$-\emph{compactly constrained disturbances}. 
Finally, we call
$$
\Gp(z_0, \UA) \mydef \sup_ {x (\cdot) \in \Xpr(z_0, \UA)} \qndx (x(\cdot)).
$$
the {\em guaranteed result} at $ z_0 \in G_0 $ for a full-memory control strategy $\UA  $ against \emph{open-loop disturbances}; and we call
$$
\Gp(z_0) \mydef \inf_{\UA \in \Uxu} \Gp(z_0, \UA).
$$
the  {\em optimal guaranteed result} at $ z_0 \in G_0 $ in $\Uxu$ against \emph{open-loop disturbances}.

Along with the full-memory control strategies, we introduce, after \cite{SubChe81e}, control quasi-strategies --- nonanticipatory transformations of disturbances into controls. The controlling player uses quasi-strategies if he/she is fully informed about the current histories and current values of the disturbance.
A {\em control quasi-strategy} is a mapping $\alpha(\cdot): {\Vpr}\mapsto {\Upr}$ satisfying the following condition:
$\alpha(v(\cdot))|_{[t_0,\tau]}= \alpha (v '(\cdot))|_{[t_0,\tau]}$ for any $ \tau\in T $, $v (\cdot),v '(\cdot) \in \Vpr$ such that $v(\cdot)|_{[t_0,\tau]} = v'(\cdot)|_{[t_0,\tau]}$.
We denote by $\Uqs$ the set of all control quasi-strategies.
For every $ z_0 \in G_0 $ and every control quasi-strategy $ \alpha(\cdot)$,  we call
$$
\Xs(z_0,\alpha(\cdot))\mydef\{x (\cdot, t_0, z_0,\alpha(v(\cdot)),v(\cdot))\mid v(\cdot)\in{\Vpr}\}
$$
the bundle
of {\em motions} originating at $ z_0 $ under $\alpha(\cdot)$.
For every
$ z_0 \in G_0 $
the value
$$
\Gq(z_0, \alpha(\cdot))=\sup_{x(\cdot)\in\Xs(z_0,\alpha(\cdot))}\qndx(x(\cdot))
$$
is called the \emph{guaranteed result} at $z_0 $ 
for a control-quasi-strategy
$ \alpha(\cdot) $
against
{\em arbitrary disturbances},
and
$$
\Gq(z_0) \mydef\myinf_{\alpha(\cdot)\in\Uqs}\Gq(z_0, \alpha(\cdot)))
$$
is called the \emph{optimal guaranteed result} at $z_0 $ \emph{in the class of the control quasi-strategies}, $\Uqs$, against {\em arbitrary disturbances}.

\begin{thm}
For every $z_0\in G_0 $
\beq\label{neq_gam}
\Gq(z_0) \leq \Gp (z_0) \leq \Gc(z_0) \leq \Gs (z_0).
\eeq
\end{thm}

\begin{rem}
From the results \cite{KraSub88e}, and \cite{SubChe81e} follows, that for every $ z_0\in G_0$ all the inequalities in (\ref{neq_gam})  turn into equalities if
\begin{equation}
\label{cond_SPSG}
\min_{u\in\PC}\max_{v\in\QC}\SP{s}{f(\tau,x,u,v)}=\max_{v\in\QC}\min_{u\in\PC}\SP{s}{f(\tau,x,u,v)}
\end{equation}
for all $(\tau,x)\in G$, $s\in\RA^n$.
In that case neither the $\Lp$-compact, nor open-loop constraints on the disturbances change the optimal guaranteed result.
\end{rem}

In this paper we do not assume (\ref{cond_SPSG})  to be satisfied 
for all $(\tau,x)\in G$, $s\in\RA^n$.
In such circumstances, some inequalities given in (\ref{neq_gam}) can be strict. 
Examples of the situations where the first and  last elements in the chain (\ref{neq_gam})
differ are well known (see \cite[Chapter VI, \S1]{SubChe81e}).
For the case where the cost functional $\qndx$  is uniformly $(L^1,\delta)$-continuous on the set of all motions of system (\ref{sys}) but is not continuous on $C(T, \RA^n)$, an example of the situation where the last inequality in (\ref{neq_gam}) is strict, was constructed in \cite{Kry91}
(where one can also find a definition of the uniform $(L^1,\delta)$-continuity).
For $\qndx$  continuous on $C(T, \RA^n)$ a similar example was given in \cite{Ser_UDSU_10e}.

Among the optimal guaranteed results 
(at a
$ z_0 \in G_0 $)
given in
 \fref{neq_gam} 
the smallest one is the optimal guaranteed result in the class of the control quasi-strategies. 
We address a question whether
the optimal guaranteed 
result 
(at 
$ z_0 $)
in the class of the full-memory control strategies
against either
open-loop disturbances,
or 
 $\Lp$-compactly constrained disturbances
 coincides with that in the class of quasi-strategies.
 If the answer is positive,
 the class of the full-memory control strategies, $\Uxu$, is 
 {\em non-improvable} against 
 a corresponding type of functional constraints on the disturbances.
In that situation, the use of any information on the past and current values of the
 actual disturbance 
does not allow the controlling player to improve the value of the 
optimal guaranteed result
at any 
$ z_0 \in G_0 $,
provided 
the disturbing player's choices are subject to
the corresponding type of functional constraints.

In \cite{Kry91} it was shown that in the case of a uniformly $ (L^1, \delta) $-continuous cost functional the one-to-one correspondence in the mapping $ v\mapsto f(t, x, u, v)$ for all $ (t,x,u) \in T \times \RA^n \times \PC $ is sufficient for the non-improvability of $\Uxu$ against $L_2$-compactly constrained disturbances. 

In \cite{Ser_DAN2013et} for the case of a cost functional continuous in $ C (T, \RA^n) $ the following sufficient condition for the non-improvability of $\Uxu$  against the $L_2$-compactly constrained disturbances was given:  for all $(t, x,u)\in G\times\PC$ we denote $\QC_{txu}$ the quotient set of the set $\QC$, generated by the equivalence relation $\mysim{txu}$: $(v_1\mysim{txu} v_2)\Leftrightarrow(f(t,x,u,v_1)=f(t,x,u,v_2))$. The condition consisted in independency of $\QC_{txu}$ on $u\in\PC$:
\beq\label{cond_Q_not_dep_u}
\QC_{txu}=\QC_{txu'}\qquad\text{for all } u,u'\in\PC,\ (t,x)\in G.
\eeq

In \cite{Ser_TRIMM2014e} for the case of a cost functional continuous in $C(T, \RA^n)$ demonstrated, that the first and second relations in \fref {neq_gam} turns into an equality, without any additional condition to be assumed.
So, it is shown that at every $ z_0 \in G_0 $ the class of the full-memory control strategies, $\Uxu$, is non-improvable against both $\Lp$-compactly constrained and open-loop disturbances.

\section {Non-improvability of full-memory control strategies}

In this section we construct a family $(\Ue)_ {\epar>0} $ of full-memory control strategies, $\Ue = (\Uepar)_{\Delta \in \Delta_T} $ $ (\varepsilon > 0) $, such that for a given $ z_0 \in G_0 $
$$
\limsup_{\varepsilon \rightarrow 0} \Gc(z_0, \Ue) \leq \Gq(z_0). 
$$
Then, in view of (\ref{neq_gam}), we get
$$
\Gc(z_0) = \Gp(z_0) = \Gq(z_0), 
$$
which implies that the full-memory control strategies are non-improvable at $ z_0 $ against both $L_2$-compactly constrained and open-loop disturbances.

The process of operation of the full-memory feedback $\Uepar= (\Uep{i}(\cdot))_{i\in\nint{0}{(\nPart-1)}} $ for a partition $\Part=\cPart$ includes on-line simulation of a motion $ y(\cdot) $ of an auxiliary copy of system (\ref{sys}), which we call the $y$-{\em model}, on every interval $ [\tau_i, \tau_{i+1}) $. In the simulation process, the control side implements the robust dynamical inversion approach (\cite{KryOsi83e,KryOsi95}). He/she identifies a 'surrogate' disturbance $ \bar{v}_i $ that mimics the affect of the actual disturbance on the system, and lets the 'surrogate' disturbance operate in the $y$-model. To identify the 'surrogate' disturbance $ \bar{v}_i $, in a small final part of the time interval $ [\tau_{i-1}, \tau_{i}) $ the controlling player implements a series of test control actions $ u^{\varepsilon}_1, \ldots, u^{\varepsilon}_{n_{\varepsilon}} $ and observes the system's reactions driven by the actual disturbance.
In the major initial part of $ [\tau_i, \tau_{i+1}) $ the controlling player implements the useful control action $ u_i $ constructed as the optimal response to the 'surrogate' disturbance for the $ y $-model, whereas the latter is driven by the useful control action $ u_{i-1} $ and 'surrogate' disturbance $\bar{v}_{i-1}$ formed previously.
The optimal response $u_i$ is found using Krasovskii's extremal shift principle (\cite{KraSub88e}); $ u_i $ shifts the $ y $-model to a target set at the maximum speed. The target set is formed in advance and  comprises the histories (up to time $ \tau_i $) of the uniform limits of the system's motions corresponding to 'approximately optimal' control quasi-strategies. The above control process ensures that the current histories of both the system's and $ y $-model's motions never abandon small neighborhoods of the current target sets, implying that at the final time, $ \vartheta $, the value of the cost functional does not exceed $\Gq(z_0) + \varphi(\varepsilon)$ for some $\varphi(\cdot)$ satisfying $\varphi(\epar)\mathop{\to}\limits_{\epar\to0}0$.
 
Now we turn to formal definitions. In the construction of the target sets we use the system's motions corresponding to 'approximately optimal' control quasi-strategies. We set
$$
\Wq(z)\mydef\bigcap_{\delta>0}\close
\bigcup_{\Gq(z, \alpha(\cdot)) \leq \Gq(z) + \delta}\Xs(z,\alpha(\cdot));
$$
here $\close X$ denotes the closure of a $X\subset C(T; \RA^n)$ in $C(T; \RA^n)$. For every $ \tau \in T $ the set of the restrictions of all the elements of $ \Wq(z) $ to $ [t_0, \tau] $, denoted by $ \Wq(z)|_{[t_0,\tau]} $, will be regarded as the {\em target set} at time $ \tau $. For every $ \tau \in T $ and every $y(\cdot)\in C([t_0,\tau],\RA^n)$ we fix a projection  $w(\cdot | \tau, y(\cdot))$ of $ y(\cdot) $ onto the target set $ \Wq(z)|_{[t_0,\tau]} $; thus,
\beq\label{pro_on_Wq}
w(\cdot |\tau, y(\cdot))\in\argmin_{w(\cdot)\in\Wq(y(t_0))|_{[t_0,\tau]}}\|w(\cdot)-y(\cdot)\|_{C([t_0,\tau],\RA^n)}.
\eeq

Fix an  $\epar\in(0,1)$. Fix an $\epar$-net $(u^\epar_j)_{j \in\nint{1}{n_\epar}}$ in $\PC$; thus, $\sup_{u\in\PC}\mymin_{j\in\nint{1}{n_\epar}}\|u-u^\epar_j\|\le\epar$. In the subsequent constructions the elements of $(u^\epar_j)_{j \in\nint{1}{n_\epar}}$ play the role of test control actions mentioned above.

Let $\Part\mydef\cPart$ be a partition of $T$. For simplicity we give the definitions for the case of partition with constant step.
Denote 
\begin{gather}
\tau_i'\mydef\tau_i-\epar \diam{\Part}, \quad i\in\nint{1}{(\nPart-1)},\label{diamin}\\
\tau'_{ij}\mydef\tau_{i}'+\frac{j(\tau_{i}-\tau_{i}')}{n_{\epar}},\quad j\in\nint{0}{n_{\epar}}, \quad i\in\nint{1}{(\nPart-1)}.\label{dop_inst}
\end{gather}
For every $x(\cdot)\in C(T;\RA^n)$, $j\in\nint{1}{n_{\epar}}$, $i\in\nint{1}{(\nPart-1)}$ let
\beq\label{ddd}
d_{ij} (x(\cdot))\mydef\frac{x(\tau'_{ij})-x(\tau'_{i(j-1)})}{\tau'_{ij }-\tau'_{i(j-1)}}.
\eeq

Define a full-memory feedback $\Uepar= (\Uep{i}(\cdot))_{i\in\nint{0}{(\nPart-1)}} $ for $\Part$ inductively. Fix some $u_*\in\PC$, $v_* \in \QC$. For every $ x_0 (\cdot) \in C([t_0,\tau_0],\RA^n)$ (recall that $\tau_0 = t_0$) we set
\begin{gather}
\bar{v}_0\mydef v_*,\quad u_0\mydef u_*,\quad  y_0(\tau_0)=z_0,\label{u_0_bar_v_0}\\
\Uep{0}(x_0(\cdot))(t)\mydef
\begin{cases}
u_0,&t\in[\tau_0,\tau_1'),\\
u^\epar_j,&t\in[\tau'_{1(j-1)},\tau'_{1j}), \ j\in\nint{1}{n_\epar}.
\end{cases}\label{Ulp0}
\end{gather}

If for some $i\in\nint{1}{(\nPart-1)}$ elements 
$$
\bar v_{i-1}=\bar v_{i-1}(x_{i-1}(\cdot))\in\QC,\quad
\Uep{(i-1)}(x_{i-1}(\cdot)) \in {\cal U}|_{[\tau_{i-1},\tau_{i}]},
$$
and
$$
y_{i-1}(\cdot)=y_{i-1}(\cdot,x_{i-1}(\cdot))\in C([t_0,\tau_{i-1}], \RA^n)
$$  
(a motion of the $y$-model on $[t_0, \tau_{i-1}]$) are defined for all $x_{i-1}(\cdot)\in C([t_0, \tau_{i-1 }], \RA^n)$, then for every $x_{i}(\cdot)\in C([t_0,\tau_i], \RA^n)$ we define $y_i(\cdot) =y_{i}(\cdot,x_{i}(\cdot))\in C([t_0,\tau_{i}], \RA^n)$ as the extension of $y_{i-1}(\cdot)$ to $[t_0, \tau_{i}]$ such that
\begin{multline}
y_{i}(\tau)=y_{i-1}(\tau_{i-1},x_i(\cdot)|_{[t_0,\tau_{i-1}]})\\
+\int_{\tau_{i-1}}^{\tau} f(t,y_{i}(t),\Uep{i-1}(x_i(\cdot)|_{[t_0,\tau_{i-1}]})(\tau_{i-1}),\\
\bar{v}_{i-1}(x_i(\cdot)|_{[t_0,\tau_{i-1}]}))dt,\quad
\tau\in[\tau_{i-1},\tau_i],\label{mov_y_def}
\end{multline}
and set
\begin{gather}
\bar{v}_{i}\in\argmin_{v\in\QC}\mymax_{j\in\nint{1}{n_\epar}}\|d_{ij}(x_i(\cdot))-f(\tau_i,x_{i}(\tau_i),u_j^\epar,v)\|,
\label{bar_v_def}\\
u_{i}\in\argmin_{u \in\PC}\SP{y_{i}(\tau_i)-w(\tau_i\mid\tau_i, y_i(\cdot))}{f(\tau_i,y_{i}(\tau_i),u,\bar{v}_{i})},\label{U_bPart_def}\\
\Uep{i}(x_i(\cdot))(t)\mydef
\begin{cases}
u_i,&t\in[\tau_i,\tau_{i+1}'),\\
u^\epar_j,&t\in[\tau'_{(i+1)(j-1)},\tau'_{(i+1)j}), j\in\nint{1}{n_\epar}.
\end{cases}
\label{Ulpi}
\end{gather}
The full-memory feedback $\Uepar$ is defined for the partition $\Part\in\Pparset$. Thus, the full-memory strategy $\Ue\mydef(\Uepar)_{\Part\in\Pparset}$ is defined.

\begin{figure}
\begin{center}
\input{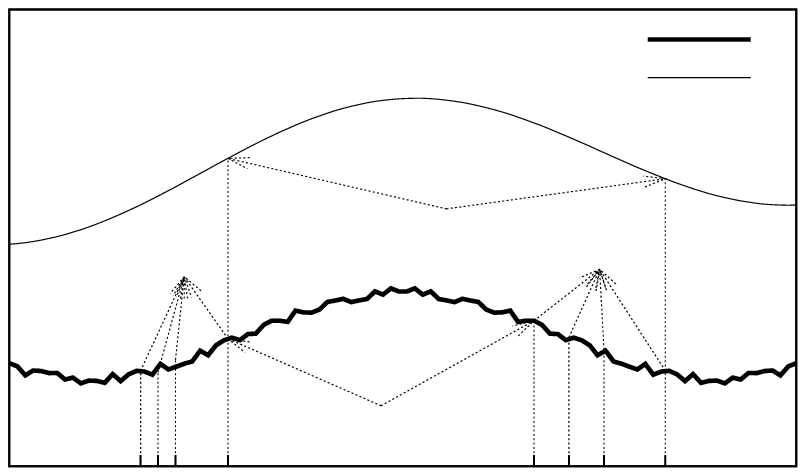}
\end{center}
\caption{The scheme of the strategy $\Ue$.}\label{pic-Ue}
\end{figure}

Illustration of the proposed control scheme is shown in Figure \ref{pic-Ue}.

As it was mentioned, the following theorem holds true.

\begin{thm}\label{teo_Q_EQ_C}
For all $z_0\in G_0$ the relations
\beq\label{Ue_approx_gamma_ps}
\limsup_{\epar\to0}\Gc(z_0,\Ue)\le\Gq(z_0),
\eeq
\beq\label{gamma_pr_eq_low_game_value}
\Gp(z_0) = \Gc(z_0) = \Gq(z_0)
\eeq
are fulfilled.
\end{thm}

\section {The constructive modifications of the optimal strategy}

In the construction of the strategy $\Ue$ there are at least two places, that may constitute essential difficulties when trying to implement numerically this control procedure.

The first is related to the calculation of projections of movements of $y$-model onto the 'target' set (see \fref{pro_on_Wq}). Conceptually the problem reduces to calculation of the gradient of 'lower' (maximin) game value in the current state of the control system. Despite the difficulty of this problem, it has long been known, extensively investigated and has viable solutions in many important cases.

The second difficulty is unlimited and fairly rapid growth of sets $(u^\epar_j)_{j\in\nint{1}{n_\epar}}$ with decreasing of parameter $\epar$. This leads to significant increase in the dimension of minimization problem \fref{bar_v_def}. Below a sufficiently broad class of systems is provided, wherein the difficulty can be overcome.

For any $(\tau, x, u, v)\in G\times\PC\times\QC$ denote the $q_{txu}(v)$ the (unique) element of $\QC_{txu}$, containing $v$.

\begin{assum}\label{con1}
There is a finite subset $\{\bar u_j \in\PC\mid j \in\nint{1} {l}\}$ such that for all $(\tau, x, u, v)\in G\times\PC\times\QC$ holds
\beq\label{bm}
\bigcap_{j\in\nint{1}{l}}q_{tx\bar u_j}(v)\subseteq q_{txu}(v).
\eeq
\end{assum}

\begin{rem}
The  assumption implies that for every $v\in\QC$ response of the system to the control $u\in\PC$ can be calculated knowing the response of the system at the same $v$ on the final set of test control actions $\{\bar u_j \in\PC \mid j \in \nint{1}{l} \}$. And so, to select an approximating value $\bar v$ (see \fref{bar_v_def}) this final set is enough. It is easy to see that this condition generalizes the condition \fref{cond_Q_not_dep_u}.
\end{rem}

We define a family of strategies $(\bar\UA_\epar)_{\epar}$ ($\bar\UA_\epar\in\Uxu$, $\epar>0$), $\bar\UA_\epar=(\bar\UB^\Part_\epar)_{\Part\in\Pparset}$, where for every $\Part \in \Pparset $ the feedback with full memory $\bar\UB^\Part_\epar$ is defined by relations \fref{diamin}--\fref{Ulpi}, wherein $n_\epar\mydef l$  and $u^\epar_j\mydef\bar u_j$, $j\in\nint{1}{n_\epar}$.

\begin{thm}\label{teo_fin_test}
Let the controlled system \fref{sys} satisfies Assumption \ref{con1}. Then, for all $ z_0 \in G_0 $ the following equalities hold
\beq\label{bUe_approx_gamma_ps}
\limsup_{\epar\to0}\Gc(z_0,\bar\UA_\epar)=\Gc(z_0).
\eeq
\end{thm}

\begin{rem}
From the construction it is clear that in case of the conditions from the theorem \ref{teo_fin_test}, in the problem of inverse dynamics \fref {bar_v_def} the data amount is fixed.
\end{rem}

The proof of Theorem \ref{teo_fin_test} in its basic steps follows the proof of Theorem \ref{teo_Q_EQ_C}.

Another control strategy $\Ul$ (see \cite[Theorem 2]{Ser_DAN2013et})  uses the value of control at the previous step as the only 'test control action'. In terms of this work that means $\tau'_i\mydef\tau_{i-1}$, $n_\epar\mydef1$, $u^\epar_1=u_{i-1}$. Thanks to condition \fref{cond_Q_not_dep_u} it was enough to identify a surrogate disturbance. 
This construction can be generalized by using instead of condition \fref{cond_Q_not_dep_u} the following assumption: 
\begin{assum}\label{con2}
There exists a closed subset $\bar\PC\subseteq\PC$ such that for all $(\tau, x, v, s)\in G\times\QC\times\RA^n$, $u,u'\in\bar\PC$ the following relations hold
\begin{gather}
\argmin_{u\in\PC}\SP{s}{f(\tau,x,u,v)}\cap\bar\PC\neq\varnothing,\nonumber\\
\QC_{txu}=\QC_{txu'}\label{sup_sub_set}.
\end{gather}
\end{assum}

Let define this modification of the strategy $\Ul\mydef(\Ulpar)_{\Part\in\Pparset}$ formally: define a full-memory feedback $\Ulpar= (\Ulp{i}(\cdot))_{i\in\nint{0}{(\nPart-1)}}$ for $\Part$ inductively. For every $ x_0 (\cdot) \in C([t_0,\tau_0],\RA^n)$ we set
\begin{gather}
\bar{v}_0\mydef v_*,\quad u_0\mydef u_*,\quad  y_0(\tau_0)=z_0,\label{u_0_bar_v_0_L}\\
\Ulp{0}(x_0(\cdot))(t)\mydef u_0,\quad t\in[\tau_0,\tau_1).
\end{gather}

\begin{figure}
\begin{center}
\input{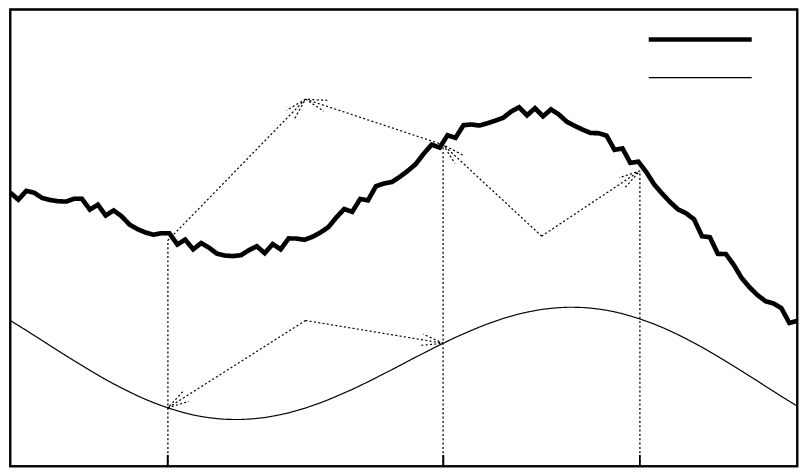}
\end{center}
\caption{The scheme of the strategy $\Ul$.}\label{pic-Ul}
\end{figure}

If for some $i\in\nint{1}{(\nPart-1)}$ elements 
$$
\bar v_{i-1}=\bar v_{i-1}(x_{i-1}(\cdot))\in\QC,\quad
\Ulp{(i-1)}(x_{i-1}(\cdot)) \in\Upr|_{[\tau_{i-1},\tau_{i}]},
$$
and
$$
y_{i-1}(\cdot)=y_{i-1}(\cdot,x_{i-1}(\cdot))\in C([t_0,\tau_{i-1}], \RA^n)
$$  
are defined for all $x_{i-1}(\cdot)\in C([t_0, \tau_{i-1 }], \RA^n)$, then for every $x_{i}(\cdot)\in C([t_0,\tau_i], \RA^n)$ we define $y_i(\cdot) =y_{i}(\cdot,x_{i}(\cdot))\in C([t_0,\tau_{i}], \RA^n)$ as the extension of $y_{i-1}(\cdot)$ to $[t_0, \tau_{i}]$ such that
\begin{multline}
y_{i}(\tau)=y_{i-1}(\tau_{i-1},x_i(\cdot)|_{[t_0,\tau_{i-1}]})\\
+\int_{\tau_{i-1}}^{\tau} f(t,y_{i}(t),\Ulp{i-1}(x_i(\cdot)|_{[t_0,\tau_{i-1}]})(\tau_{i-1}),\\
\bar{v}_{i-1}(x_i(\cdot)|_{[t_0,\tau_{i-1}]}))dt,\quad
\tau\in[\tau_{i-1},\tau_i],\label{mov_y_def_L}
\end{multline}
and set
\begin{multline}\label{bar_v_def_L}
\bar{v}_{i}\in\argmin_{v\in\QC}\Bigl\|\frac{x_i(\tau_i)-x_i(\tau_{i-1})}{\tau_i-\tau_{i-1}}\\
-f(\tau_i,x_{i}(\tau_i),\Ulp{i-1}(x_i(\cdot))(\tau_{i-1}),v)\Bigr\|,
\end{multline}
\begin{gather}
u_{i}\in\argmin_{u \in\bar\PC}\SP{y_{i}(\tau_i)-w(\tau_i\mid\tau_i, y_i(\cdot))}{f(\tau_i,y_{i}(\tau_i),u,\bar{v}_{i})},\label{U_bPart_def_L}\\
\Ulp{i}(x_i(\cdot))(t)\mydef u_i,\qquad t\in[\tau_i,\tau_{i+1}).\label{Ulpi_L}
\end{gather}
The full-memory feedback $\Ulpar$ is defined for the partition $\Part\in\Pparset$. So, the full-memory strategy $\Ul\mydef(\Ulpar)_{\Part\in\Pparset}$ is defined. The  scheme of the control strategy is shown in Figure \ref{pic-Ul}.

\begin{thm}\label{teo_one_test}
Let the controlled system \fref{sys} satisfies Assumption \ref{con2}. Then, the equalities 
\beq\label{Ul_get_gamma_qs}
\Gc(z_0,\Ul)=\Gc(z_0)
\eeq
are fulfilled for all $z_0\in G_0$.
\end{thm}

\section{Example}

Let system \fref{sys} have the form
\beq
\label{sys_examp_2x2}
\begin{cases}
\dot x_1(\tau) = u_1(\tau) v_1(\tau),\quad\tau\in T\mydef[0,1],\\
\dot x_2(\tau) =\max\{0, x_1(\tau)\}u_2(\tau) v_2(\tau),\\
(x_1(0),x_2(0))=(0,0),\quad G_0 = \{(0,0)\},
\end{cases}
\eeq
$$
u_1(\tau), u_2(\tau)\in  [-1,1],\quad v_1(\tau),v_2(\tau)\in \{-1,1\},
$$
and the cost functional be given by $\gamma(x(\cdot))\mydef x_2(1)$ $(x(\cdot) = (x_1(\cdot), x_2(\cdot)) \in C(T,\RA^2))$. With an appropriate choice of
$G$, system \fref{sys_examp_2x2} satisfies all the assumptions imposed earlier on system \fref{sys}; therefore, Theorem \ref{teo_Q_EQ_C} holds, implying the full-memory control strategies are non-improvable against both the $L_2$-compactly constrained and open-loop disturbances. On the other hand, system \fref{sys} does not satisfy the conditions sufficient for the non-improvability of the full-memory control strategies against the $L_2$-compactly constrained disturbances, which are given in \cite{Kry91}  (Theorem 9.1) and in \cite{Ser_DAN2013et} (Theorem 2). By use of relations \fref{gamma_pr_eq_low_game_value} one can find that $\Gc((0,0))=\Gp((0,0))=\Gq((0,0))=-0.5$. 

It is clear that the set $\bar\PC\mydef\{-1,1\}\subset\PC$ satisfies to both assumptions \ref{con1}, \ref{con2}. So, we can use the construction of optimal strategy $\UA_*=(\UB_*^\Part)_{\Part\in\Pparset}$ given in \fref{u_0_bar_v_0_L}--\fref{Ulpi_L}: by using the monotonicity of the quality index $\qndx$, we get the feedback with full memory $\Ulpar=(\Ulp{i})_{i\in\nint{0}{\nPart-1}}$ for the partition $\Part\mydef\cPart$, the motion $(x(\cdot)\in C([t_0,\tau_i],\RA^2))$ and $i\in\nint{1}{(\nPart-1)}$:
\beqnt
\Ulp{i}(x(\cdot))\in\LL{(}\argmax\limits_{u_1\in\bar\PC}\{u_1\cdot\frac{x_1(\tau_{i})-x_1(\tau_{i-1})}{u_1(\tau_{i-1})}\}\atop\argmin\limits_{u_2\in\bar\PC}\{u_2\cdot\frac{x_2(\tau_{i})-x_2(\tau_{i-1})}{u_2(\tau_{i-1})}\}\RR{)},
\eeq
$$
(u_1(\tau_{i-1}),u_2(\tau_{i-1}))\mydef\Ulp{(i-1)}(x(\cdot)|_{[t_0,\tau_{i-1}]}).
$$

\begin{ack}
This work was supported by the Russian Foundation for Basic Research (project no.  12-01-00290),  by the Program for Fundamental Research of Presidium of the Russian Academy of Sciences ``Dynamic Systems and Control Theory'', by the Ural Branch of the Russian Academy of Sciences (project no. 12-$\Pi$-1-1002).
\end{ack}


\end{document}